\documentclass[12pt]{amsart}
\usepackage{cases}
\usepackage{txfonts}
\usepackage{latexsym}
\usepackage{amsthm}
\usepackage{mathrsfs}
\usepackage{amssymb, amsmath}

\def\opn#1#2{\def#1{\operatorname{#2}}} % to make operators
\opn\chara{char} \opn\length{\ell}
\opn\projdim{proj\,dim} \opn\injdim{inj\,dim} \opn\rank{rank}
\opn\depth{depth} \opn\grade{grade} \opn\height{height}
\opn\embdim{emb\,dim} \opn\codim{codim}

\opn\Tr{Tr} \opn\bigrank{big\,rank}
\opn\superheight{superheight}\opn\lcm{lcm}
\opn\trdeg{tr\,deg}%
\opn\reg{reg} \opn\lreg{lreg}
\opn\div{div} \opn\Div{Div} \opn\cl{cl} \opn\Cl{Cl}
\opn\Spec{Spec} \opn\Supp{Supp} \opn\supp{supp} \opn\Sing{Sing}
\opn\Ass{Ass}
\opn\Ann{Ann} \opn\Rad{Rad} \opn\Soc{Soc}
\opn\Ker{Ker} \opn\Coker{Coker} \opn\Im{Im} \opn\Hom{Hom}
\opn\Tor{Tor} \opn\Ext{Ext} \opn\End{End} \opn\Aut{Aut} \opn\id{id}

\opn\nat{nat}
\opn\pff{pf}%   \pf exists already
\opn\Pf{Pf} \opn\GL{GL} \opn\SL{SL} \opn\mod{mod} \opn\ord{ord}
\opn\aff{aff} \opn\con{conv} \opn\relint{relint} \opn\st{st}
\opn\lk{lk} \opn\cn{cn} \opn\core{core} \opn\vol{vol}
\opn\gr{gr}

\def\pot#1#2{#1[\kern-0.28ex[#2]\kern-0.28ex]}

\opn\dirlim{\underrightarrow{\lim}}
\opn\invlim{\underleftarrow{\lim}}
\newtheorem{theorem}{Theorem}

\newtheorem{lemma}[theorem]{Lemma}

\def\qed{\ifhmode\textqed\fi
   \ifmmode\ifinner\quad\qedsymbol\else\dispqed\fi\fi}
\def\textqed{\unskip\nobreak\penalty50
    \hskip2em\hbox{}\nobreak\hfil\qedsymbol
    \parfillskip=0pt \finalhyphendemerits=0}
\def\dispqed{\rlap{\qquad\qedsymbol}}

\opn\ini{in} \opn\inm{inm} \opn\Sym{Sym}

\begin{document}
\title{distinct zeros and simple zeros of Dirichlet $L$-functions}
\author{XiaoSheng Wu}
\date{}
\address {School of Mathematical Sciences\\ University of Science and Technology of China\\ Hefei 230026, P. R.
China.}
\email {xswu@amss.ac.cn}
\subjclass[2010]{11M26, 11M06 }
\keywords{simple zeros; distinct zeros; Dirichlet $L$-function.}

\begin{abstract} In this paper, we study the number of additional zeros of Dirichlet $L$-function caused by multiplicity by using Asymptotic Large Sieve. Then in asymptotic terms we prove that there are more than 80.13\% of zeros of the family of Dirichlet $L$-functions are distinct and more than 60.261\% of zeros of the family of Dirichlet $L$-functions are simple. In addition, assuming the Generalized Riemann Hypothesis, we improve these proportions to 83.216\% and 66.433\%.

\end{abstract}
\maketitle

\section{Introduction}
Let $L(s,\chi)$ be a Dirichlet $L$-function, where $s=\sigma+it$, $\chi$ (mod $q$) is a character. It is defined for $\sigma>1$ by
\begin{align}
   L(s,\chi)=\sum_{n=1}^\infty\chi(n)n^{-s}.
\end{align}
Let $\chi$ (mod $q$) be a primitive character. The total number of zeros $\rho=\beta+i\gamma$ of $L(s,\chi)$ with $0<\beta<1$ and $|\gamma|\le T$, say $N(T,\chi)$, is known asymptotically very precisely
\begin{align}
   N(T,\chi)=\frac{T}{\pi}\log\frac{qT}{2\pi e}+O(\log qT),\ \ \ \ \ T\ge3.
\end{align}

We define the number of distinct zeros and the number of simple zeros of a single $L(s,\chi)$ as follows
\begin{align}
   &N_d(T,\chi)=|\{\rho=\beta+i\gamma: -T<\gamma\le T,~L(\rho,\chi)=0\}|,\notag\\
   &N_s(T,\chi)=|\{\rho=\beta+i\gamma: -T<\gamma\le T,~L(\rho,\chi)=0, L'(\rho,\chi)\neq0\}|.\notag
\end{align}
It is believed that  $N_d(T,\chi)=N_s(T,\chi)=N(T,\chi)$, which means that all zeros of Dirichlet $L$-function are simple. This is known as the Simple Zero Conjecture.

As a special case, for the Riemann zeta-function, there is a long history of the study of the Simple Zero Conjecture, and the topic has been studied by a lot of papers (to see \cite{And, Con-More, CG, Far, Hea} for example).

In 1995, Farmer \cite{Far} proved that at least 63.952\% of zeros of the Riemann zeta-function are distinct. Farmer proved this by using a combination method which was  based on proportions of simple zeros of $\xi^{(n)}(s,1)$.

It is known that more than 40.58\% of zeros of the Riemann zeta-function are simple, which was proved in 2011 by Bui, Conrey and Young's work \cite{BCY}. This work was based on Levinson's method (to see \cite{Con-Zero, Con-More, CG, Lev-More}) with a more general mollifier. If assumed the Riemann Hypothesis, Cheer and Goldston proved in 1993 that more than 67.275\% of zeros of the Riemann zeta-function are simple.

In 1998, under the assumption of the Riemann Hypothesis and the Generalized Lindel\"{o}f Hypothesis, Conrey, Ghosh and Gonek proved in \cite{CGG} that more than 84.56\% of zeros of the Riemann zeta-function are distinct and more than 70.37\% of zeros of the Riemann zeta-function are simple. Their starting point is the observation that
\begin{align}
\bigg|\sum_{-T\le\gamma\le T}M(\frac12+i\gamma)\zeta'(\frac12+i\gamma)\bigg|^2\le N_s(T,1)\sum_{-T\le\gamma\le T}\big|M(\frac12+i\gamma)\zeta'(\frac12+i\gamma)|^2\notag
\end{align}
by Cauchy's inequality. The function $M(s)$ is taken to be a mollifier
\begin{align}
   M(s)=\sum_{n\le y}\mu(n)P\bigg(\frac{\log y/n}{\log y}\bigg)n^{-s},\notag
\end{align}
where $y=T^{1/2-\epsilon}$ and $P(x)$ is a suitable polynomial.
%By a combination method based on proportions of simple zeros of $\xi^{(n)}(s,1)$, Farmer \cite{Far} proved in 1995 that at least 63.9\% of the zeros of the Riemann zeta-function are distinct. Based on Levinson's method (to see \cite{Con-Zero, Con-More, CG, Lev-More}), Bui, Conrey and Young \cite{BCY} proved that more than 40.58\% of zeros of the Riemann zeta-function are simple, and this is the best unconditional result on proportion of simple zeros now. Assuming the Riemann Hypothesis, the best result now is that more than 67.275\% of zeros of the Riemann zeta-function are simple. This result was proved in \cite{CG} by Cheer and Goldston, and they obtained this by improving Montgomery and Taylor's works \cite{Mon, MT}. If the Riemann Hypothesis and the Generalized Lindel\"{o}f Hypothesis hold, Conrey, Ghosh and Gonek proved in 1998 that more than 70.37\% of zeros of the Riemann zeta-function are simple and more than 84.56\% of zeros of the Riemann zeta-function are distinct. However, for general Dirichlet $L$-function with $\chi\neq1$, there is hardly any similar result.

For general Dirichlet $L$-function with $\chi\neq1$, there was few similar result. In Conrey, Iwaniec and Soundararajan's recent work \cite{CIS-Crit}, they proved that at least 58.65\% of zeros of the family of Dirichlet $L$-functions are on the critical line (1/2-line) and simple by using Asymptotic Large Sieve (to see \cite{CIS-Asym}) and Levinson's method.

In this work, we introduce a new way to  work on the number of distinct zeros and simple zeros of the family of Dirichlet $L$-functions. Here we do not care about the relationship among zeros of different Dirichlet $L$-functions. In other words, if $\rho_1$ is a zero of $L(s,\chi_1)$ and $\rho_2$ is a zero of $L(s,\chi_2)$ with $\chi_1\neq\chi_2$, we say that $\rho_1,~\rho_2$ are two different zeros of the family of Dirichlet $L$-functions even though $\rho_1=\rho_2$ numerically.

Let $\Psi(x)$ be a non-negative function defined as in \cite{CIS-Crit}, which is smooth, compactly supported on $\mathbb{R}^+$. Put
\begin{align}
\label{t1.1}
   N(T,Q)=\sum_{q\le Q}\frac{\Psi(q/Q)}{\phi(q)}\mathop{{\sum}^*}\limits_{\chi(\text{mod} q)}N(T,\chi),
\end{align}
where $Q\ge3$ and $T\ge3$. Here the superscript $*$ restricts the summation to the primitive characters.
Accordingly, we also define
\begin{align}
   N_d(T,Q)=\sum_{q\le Q}\frac{\Psi(q/Q)}{\phi(q)}\mathop{{\sum}^*}\limits_{\chi(\text{mod} q)}N_d(T,\chi),\notag
\end{align}
and
\begin{align}
   N_s(T,Q)=\sum_{q\le Q}\frac{\Psi(q/Q)}{\phi(q)}\mathop{{\sum}^*}\limits_{\chi(\text{mod} q)}N_s(T,\chi)\notag
\end{align}
for $Q\ge3$ and $T\ge3$.

In the following we will first study the number of additional zeros of Dirichlet $L$-function caused by multiplicity. Here the number of additional zeros caused by multiplicity is the number of zeros with each zero counted according to multiplicity minus one. Then we will obtain that there are more than 80.13\% of zeros of the family of Dirichlet $L$-functions are distinct and more than 60.261\% of zeros of the family of Dirichlet $L$-functions are simple. In addition, if assume that the Generalized Riemann Hypothesis holds, we can obtain better results both on distinct zeros and simple zeros. For convenience, we state our main results as follows.
\begin{theorem}
\label{thm1}
For $Q$ and $T$ with $(\log Q)^6\le T\le(\log Q)^A$ we have
\begin{align}
    N_d(T,Q)\ge0.8013 N(T,Q),\ \ \ \ \ \ \ \ \ \ \ \ \ \ N_s(T,Q)\ge0.60261 N(T,Q),
\end{align}
where $A\ge6$ is any constant, provided $Q$ is sufficiently large in terms of $A$.
\end{theorem}
If the Riemann Hypothesis for the family of these Dirichlet $L$-functions holds, we obtain the following result.
\begin{theorem}
\label{thm2}
Assume the Generalized Riemann Hypothesis. For $Q$ and $T$ with $(\log Q)^6\le T\le(\log Q)^A$ we have
\begin{align}
     N_d(T,Q)\ge0.83216 N(T,Q),\ \ \ \ \ \ \ \ \ \ \ \ \ \ N_s(T,Q)\ge0.66433 N(T,Q),
\end{align}
where $A\ge6$ is any constant, provided $Q$ is sufficiently large in terms of $A$.
\end{theorem}
%{\bf Acknowledgements.} I am very grateful to Professor Banghe Li and Professor Shaoji Feng,  my thesis advisors, for their encouragement on my work. I also want to thank the lecture ``Ergodic Prime Number Theorems Seminar 2012''.

\section{sketch of the proof}
To study the number of distinct zeros and simple zeros of Dirichlet $L$-function, we firstly focus on the number of the additional zeros caused by multiplicity. We note that if $\rho$, a zero of $L(s,\chi)$, is a non-simple zero, it must be a zero of
\begin{align}
\label{2.1}
   G(s,\chi)=L(s,\chi)\psi_1(s,\chi)+\lambda L'(s,\chi)\psi_2(s,\chi)
\end{align}
with multiplicity reduced by at most one. Here $\psi_1(s,\chi)$, $\psi_2(s,\chi)$ can be any analytic functions and $\lambda$ can be any constant. Hence the number of additional zeros of $L(s,\chi)$ caused by multiplicity in any region is not more than the number of zeros of $G(s,\chi)$ in the same region.

It is therefore important to find a such $G(s,\chi)$ that has less zeros in the critical area ($0<\text{Re}(s)<1$), but this is very hard. A more feasible way is to partition the critical area into some sub areas, and for each sub area we may find a $G(s,\chi)$ that has less zeros in the area. Then we have the number of additional zeros of $L(s,\chi)$ caused by multiplicity in the critical  area are not more than the sum of zeros of these $G(s,\chi)$ in the corresponding sub area.

In this paper, we actually partition the critical area into the left part ($\text{Re}(s)<1/2$) and the right part ($\text{Re}(s)\ge1/2$). The partition used here may be not good enough, it is possible for one to find a more useful partition and the work will be worthy.

Firstly, let us consider the left side ($\text{Re}(s)<1/2$). Actually, the estimate of additional zeros in the left side is only required in the proof of Theorem \ref{thm1}. In this side, we choose $G(s,\chi)=\xi'(s,\chi)$, where
\begin{align}
   \xi(s,\chi)=H(s,\chi)L(s,\chi)\notag
\end{align}
with
\begin{align}
\label{s2.2}
   H(s,\chi)=\frac12s(s-1)\bigg(\frac{q}{\pi}\bigg)^{\frac{s}{2}}\Gamma\bigg(\frac{s+(1-\chi(-1))/2}{2}\bigg).
\end{align}

About the number of zeros of $\xi^{(j)}(s,\chi)$ for primitive character $\chi$, we may have the following lemma, which can be proved similarly as Lemma 2 in \cite{Con-Zero}.
\begin{lemma}
\label{lems1}
Let $\chi$ be a primitive character. For any integer $j\ge0$, all zeros of $\xi^{(j)}(s,\chi)$ satisfy $0<\sigma<1$. Let $N_{\xi^{(j)}}(T,\chi)$ denote the number of zeros of $\xi^{(j)}(s,\chi)$ with $-T\le t\le T$, then
\begin{align}
   N_{\xi^{(j)}}(T,\chi)=\frac{T}{\pi}\log\frac{qT}{2\pi e}+O_j(\log qT).\notag
\end{align}
\end{lemma}
The functional equation for $L(s,\chi)$ says
\begin{align}
\label{s2.3}
   h(\chi)\xi(s,\chi)=\bar h(\chi)\xi(1-s,\bar\chi),
\end{align}
where
\begin{align}
   \frac{\bar h(\chi)}{h(\chi)}=\frac{\tau(\chi)}{i^{(1-\chi(-1))}q^{1/2}}\notag
\end{align}
with $|h(\chi)|=1$, $\tau(\chi)=\sum_{n=1}^q\chi(n)\exp(2\pi in/q)$.

From (\ref{s2.3}) we may find that if $\rho$ is a zero of $\xi'(s,\chi)$, then $1-\rho$ is a zero of $\xi'(s,\overline{\chi})$. From this and Lemma \ref{lems1} with $j=1$, we can see that $N_{\xi',l}(T,\chi)$, the number of zeros of $\xi'(s,\chi)$ , and $N_{\xi',l}(T,\overline{\chi})$, the number of zeros of $\xi'(s,\overline{\chi})$ in the left side, satisfy
\begin{align}
\label{t2.1}
   N_{\xi',l}(T,\chi)+N_{\xi',l}(T,\overline{\chi})=N(T,\chi)-N_{\xi',c}(T,\chi).
\end{align}
Here $N_{\xi',c}(T,\chi)$ denotes the number of zeros of $\xi'(s,\chi)$ on the critical line, which satisfies $N_{\xi',c}(T,\chi)=N_{\xi',c}(T,\overline{\chi})$. In particular for $\chi$ is a real character, it follows
\begin{align}
   N_{\xi',l}(T,\chi)=\frac12\Big(N(T,\chi)-N_{\xi',c}(T,\chi)\Big).\notag
\end{align}
Hence it is important for us to obtain a lower bound for the number of zeros of $\xi'(s,\chi)$ on the critical line.

We now come to the right side. Let $\mathbf{q}=q/\pi$ as in \cite{CIS-Crit} and $R>0$ be a constant to be specified later. Choosing $\lambda=1/(r\log\mathbf{q})$, let $G(s,\chi)$ be defined by (\ref{2.1}) with
\begin{align}
\label{2.30}
   &\psi_1(s,\chi)=\sum_{n\le X}\frac{\mu(n)\chi(n)}{n^{s+R/\log\mathbf{q}}}P_1\bigg(\frac{\log X/n}{\log X}\bigg),\notag\\
   &\psi_2(s,\chi)=\sum_{n\le X}\frac{\mu(n)\chi(n)}{n^{s+R/\log\mathbf{q}}}P_2\bigg(\frac{\log X/n}{\log X}\bigg),
\end{align}
where $\mu$ is the M\"{o}bius function, $X=\mathbf{q}^\theta$ with $0<\theta<1$. Here $P_1,~P_2$ are polynomials with $P_1(0)=P_2(0)=0$, $P_1(1)=P_2(1)=1$, which will be specified later. It is worth remarking that this choice of $\psi_i$ is not the best, some other possible choices of $\psi_i$ can be found in \cite{And, BCY, Feng, Lev-Dedu}.

Let $\mathcal{D}$ be the closed rectangle with vertices $1/2-iT$, $3-iT$, $1/2+iT$, $3+iT$. Let $N_G(\mathcal{D},\chi)$ denote the number of zeros of $G(s,\chi)$ in $\mathcal{D}$, including zeros on the left boundary. It is obvious that zeros of $G(s,\chi)$ in the right side of the critical area are not more than zeros of it in $\mathcal{D}$, and so the the number of additional zeros of $L(s,\chi)$ caused by multiplicity in the right side is not more than $N_G(\mathcal{D},\chi)$.

From the above discussion, it follows that the number of additional zeros of $L(s,\chi)$ caused by multiplicity in the left side is not more than the number of zeros of $\xi'(s,\chi)$ in the left side, and in the right side is not more than $N_G(\mathcal{D},\chi)$. So, by (\ref{t2.1}) we may have the following formula about the number of distinct zeros of $L(s,\chi)$
\begin{align}
\label{++1}
   N_d(T,\chi)+N_d(T,\overline{\chi})\ge N(T,\chi)+N_{\xi',c}(T,\chi)-N_G(\mathcal{D},\chi)-N_G(\mathcal{D},\overline{\chi}).
\end{align}
In addition, the number of all non-simple zeros counted according to multiplicity are not more than $2N_{\xi',l}(T,\chi)+2N_G(\mathcal{D},\chi)$ since each non-simple zero has multiplicity at least 2. Thus for simple zeros we have
\begin{align}
\label{++2}
   N_s(T,\chi)+N_s(T,\overline{\chi})&\ge 2N(T,\chi)-2N_{\xi',l}(T,\chi)-2N_G(\mathcal{D},\chi)-2N_{\xi',l}(T,\overline{\chi})-2N_G(\mathcal{D},\overline{\chi})\notag\\
   &=2N_{\xi',c}(T,\chi)-2N_G(\mathcal{D},\chi)-2N_G(\mathcal{D},\overline{\chi}).
\end{align}

If assume the Riemann Hypothesis for this $L(s,\chi)$, we do not need to care about the left side since there is no zero in this side. Then we have
\begin{align}
   N_d(T,\chi)\ge N(T,\chi)-N_G(\mathcal{D},\chi).\notag
\end{align}
Similarly, for simple zero we have
\begin{align}
   N_s(T,\chi)\ge N(T,\chi)-2N_G(\mathcal{D},\chi).\notag
\end{align}

As we have seen, it is important therefore to obtain a lower bound for $N_{\xi',c}(T,\chi)$ and an upper bound for $N_G(\mathcal{D},\chi)$. Unfortunately, at the present state of technology we are unable to obtain a useful evaluation of them for individual characters $\chi~(\text{mod} q)$. However, by averaging over the conductors $q$ and the primitive characters $\chi$, we are able to get a useful evaluation of them by using Asymptotic Large Sieve.

As the definition of $N(T,Q)$ in (\ref{t1.1}) we define
\begin{align}
   N_{G}(\mathcal{D},Q)=\sum_{q\le Q}\frac{\Psi(q/Q)}{\phi(q)}\mathop{{\sum}^*}\limits_{\chi(\text{mod} q)}N_{G}(\mathcal{D},\chi),\notag
\end{align}
and
\begin{align}
   N_{\xi',c}(T,Q)=\sum_{q\le Q}\frac{\Psi(q/Q)}{\phi(q)}\mathop{{\sum}^*}\limits_{\chi(\text{mod} q)}N_{\xi',c}(T,\chi),\notag
\end{align}
for $Q\ge3$ and $T\ge3$. Then we have from (\ref{++1}) and (\ref{++2}) that
\begin{align}
\label{2.2}
   N_d(T,Q)\ge \frac12N(T,Q)+\frac12N_{\xi',c}(T,Q)-N_G(\mathcal{D},Q),
\end{align}
\begin{align}
\label{2.21}
   N_s(T,Q)\ge N_{\xi',c}(T,Q)-2N_G(\mathcal{D},Q).
\end{align}
If assume Generalized Riemann Hypothesis, we have
\begin{align}
\label{2.3}
   N_d(T,Q)\ge N(T,Q)-N_G(\mathcal{D},Q),
\end{align}
\begin{align}
\label{2.31}
   N_s(T,Q)\ge N(T,Q)-2N_G(\mathcal{D},Q).
\end{align}

In the following,  we will first introduce an asymptotic large sieve result in section \ref{sec3}. Then we will obtain a upper bound for $N_G(\mathcal{D},Q)$ in section \ref{sec4} and a lower bound for $N_{\xi',c}(T,Q)$ in section \ref{sec5}.

\section{asymptotic large sieve result}
\label{sec3}
To evaluate $N_G(\mathcal{D},Q)$ and $N_{\xi',c}(T,Q)$, we need to evaluate the following sum
\begin{align}
  &\sum_{q\le Q}\frac{\Psi(q/Q)\log\mathbf{q}}{\phi(q)}\notag\\
  & \ \ \ \ \ \ \ \ \ \ \ \ \ \ \times\mathop{{\sum}^*}\limits_{\chi(\text{mod} q)}\int_{-T}^TL(\sigma_0+\alpha+it,\chi)L(\sigma_0+\beta-it,\bar{\chi}) \psi_{i_1}(\sigma_0+it,\chi)\psi_{i_2}(\sigma_0-it,\bar{\chi})dt\notag
\end{align}
for $\alpha,~\beta\ll1/\log Q$.
We will evaluate this sum from the following lemma, which was implied in the work of Conrey, Iwaniec and Soundararajan \cite{CIS-Crit}.
\begin{lemma}
Suppose that $\psi_1(s,\chi)$, $\psi_2(s,\chi)$ are defined by (\ref{2.30}). For sufficiently large $Q$, $T$ with $\log^6Q\le T\le \log^AQ$ and constant $A\ge 6$, let $X=\mathbf{q}^\theta$ with $0<\theta<1$, $\sigma_0=1/2-R/\log \mathbf{q}$, $\alpha,\beta\ll 1/\log Q$ and $h_1=h/(h,k),~k_1=k/(h,k)$.
Then, for $i_1,i_2=1,2$,
\begin{align}
\label{3.1}
  &\sum_{q\le Q}\frac{\Psi(q/Q)\log\mathbf{q}}{\phi(q)}\notag\\
  & \ \ \ \ \ \ \ \ \ \ \ \ \ \ \times\mathop{{\sum}^*}\limits_{\chi(\text{mod} q)}\int_{-T}^TL(\sigma_0+\alpha+it,\chi)L(\sigma_0+\beta-it,\bar{\chi}) \psi_{i_1}(\sigma_0+it,\chi)\psi_{i_2}(\sigma_0-it,\bar{\chi})dt\notag\\
  &\sim 2T\bigg(\sum_{q\le Q}\Psi(\frac qQ)\log\mathbf{q}\frac{\phi^*(q)}{\phi(q)}\bigg)\sum_{h\le X}\sum_{k\le X}\frac{\mu(h)\mu(k)}{(hh_1kk_1)^{1/2}}F(\alpha,\beta,h_1,k_1)P_{i_1}\bigg(\frac{\log X/h}{\log X}\bigg)P_{i_2}\bigg(\frac{\log X/k}{\log X}\bigg),
\end{align}
where
\begin{align}
\label{3.2}
   F(\alpha,\beta,h_1,k_1)=h_1^{1/2-\sigma_0-\beta}k_1^{1/2-\sigma_0-\alpha}\zeta(2\sigma_0+\alpha+\beta)+ \frac{h_1^{\alpha+\sigma_0-1/2}k_1^{\beta+\sigma_0-1/2}}{\mathbf{q}^{2\sigma_0+\alpha+\beta-1}} \zeta(2-2\sigma_0-\alpha-\beta),
\end{align}
and $\phi^*(q)$ denotes the number of primitive characters (mod $q$) as in \cite{CIS-Crit}.
\end{lemma}

We now evaluate the sum over $h$ and $k$ in the right side of (\ref{3.1}). Making the following variable changes
\begin{align}
   a=(\sigma_0+\alpha-1/2)\log\mathbf{q}\ll1,\notag
\end{align}
\begin{align}
   b=(\sigma_0+\beta-1/2)\log\mathbf{q}\ll1,\notag
\end{align}
we find that the formula (\ref{3.2}) becomes to
\begin{align}
   f(a,b,h_1,k_1)=h_1^{-b/\log\mathbf{q}}k_1^{-a/\log\mathbf{q}}\zeta\bigg(1+\frac{a+b}{\log\mathbf{q}}\bigg)+ \frac{h_1^{a/\log\mathbf{q}}k_1^{b/\log\mathbf{q}}}{e^{a+b}} \zeta\bigg(1-\frac{a+b}{\log\mathbf{q}}\bigg).\notag
\end{align}
We also approximate $\zeta(s)$ near 1 by $(s-1)^{-1}$ getting the following asymptotic values
\begin{align}
   \zeta\bigg(1+\frac{a+b}{\log\mathbf{q}}\bigg)\sim\frac{\log\mathbf{q}}{a+b},\ \ \ \ \ \ \ \ \ \ \zeta\bigg(1-\frac{a+b}{\log\mathbf{q}}\bigg)\sim-\frac{\log\mathbf{q}}{a+b}.\notag
\end{align}
Then we get
\begin{align}
\label{3.3}
   f(a,b,h_1,k_1)\sim\frac{\log\mathbf{q}}{a+b}\bigg(h_1^{-b/\log\mathbf{q}}k_1^{-a/\log\mathbf{q}}- e^{-a-b}h_1^{a/\log\mathbf{q}}k_1^{b/\log\mathbf{q}}\bigg).
\end{align}
Substituting (\ref{3.3}) into the right side of (\ref{3.1}) we have
\begin{align}
\label{3.41}
   &\sum_{h\le X}\sum_{k\le X}\frac{\mu(h)\mu(k)}{(hh_1kk_1)^{1/2}}F(\alpha,\beta,h_1,k_1)P_{i_1}\bigg(\frac{\log X/h}{\log X}\bigg)P_{i_2}\bigg(\frac{\log X/k}{\log X}\bigg)\notag\\
   &\sim\frac{\log\mathbf{q}}{a+b}\Bigg\{\sum_{h\le X}\sum_{k\le X}\frac{\mu(h)\mu(k)}{h^{1/2}k^{1/2}h_1^{1/2+b/\log\mathbf{q}}k_1^{1/2+a/\log\mathbf{q}}}P_{i_1}\bigg(\frac{\log X/h}{\log X}\bigg)P_{i_2}\bigg(\frac{\log X/k}{\log X}\bigg)\notag\\
   &-e^{-a-b}\sum_{h\le X}\sum_{k\le X}\frac{\mu(h)\mu(k)}{h^{1/2}k^{1/2}h_1^{1/2-a/\log\mathbf{q}}k_1^{1/2-b/\log\mathbf{q}}}P_{i_1}\bigg(\frac{\log X/h}{\log X}\bigg)P_{i_2}\bigg(\frac{\log X/k}{\log X}\bigg)\Bigg\}.
\end{align}
To evaluate the sum over $h$ and $k$ within the brackets we appeal to Lemma 1 of \cite{Con-More}, which says that
\begin{align}
\label{3.51}
   &\sum_{h\le X}\sum_{k\le X}\frac{\mu(h)\mu(k)}{h^{1/2}k^{1/2}h_1^{1/2+\omega_1}k_1^{1/2+\omega_2}}P_{i_1}\bigg(\frac{\log X/h}{\log X}\bigg)P_{i_2}\bigg(\frac{\log X/k}{\log X}\bigg)\notag\\
   &\sim\frac1{\log X}\int_0^1\Big(P_{i_1}'(t)+\omega_1\log XP_{i_1}(t)\Big)\Big(P_{i_2}'(t)+\omega_2\log XP_{i_2}(t)\Big)dt
\end{align}
holds uniformly in complex numbers $\omega_1,\omega_2\ll(\log X)^{-1}$.
Recall that $X=\mathbf{q}^\theta$. Substituting (\ref{3.51}) into (\ref{3.41}) we get
\begin{align}
\label{3.4}
   &\sum_{h\le X}\sum_{k\le X}\frac{\mu(h)\mu(k)}{(hh_1kk_1)^{1/2}}F(\alpha,\beta,h_1,k_1)P_{i_1}\bigg(\frac{\log X/h}{\log X}\bigg)P_{i_2}\bigg(\frac{\log X/k}{\log X}\bigg)\notag\\
   &\sim\frac{1}{(a+b)\theta}\bigg(g_{i_1,i_2}(b,a)-e^{-a-b}g_{i_1,i_2}(-a,-b)\bigg)
\end{align}
with
\begin{align}
\label{3.5}
   g_{i_1,i_2}(a,b)=\int_0^1\Big(P_{i_1}'(t)+a\theta P_{i_1}(t)\Big)\Big(P_{i_2}'(t)+b\theta P_{i_2}(t)\Big)dt.
\end{align}
Substituting (\ref{3.4}) into (\ref{3.1}) we have
\begin{align}
\label{3.6}
   &\sum_{q\le Q}\frac{\Psi(q/Q)\log\mathbf{q}}{\phi(q)}\notag\\
  & \ \ \ \ \ \ \ \times\mathop{{\sum}^*}\limits_{\chi(\text{mod} q)}\int_{-T}^TL(\sigma_0+\alpha+it,\chi)L(\sigma_0+\beta-it,\bar{\chi}) \psi_{i_1}(\sigma_0+it,\chi)\psi_{i_2}(\sigma_0-it,\bar{\chi})dt\notag\\
   &\sim \frac{2T}{(a+b)\theta} \bigg(\sum_{q\le Q}\Psi(\frac qQ)\log\mathbf{q}\frac{\phi^*(q)}{\phi(q)}\bigg) \bigg(g_{i_1,i_2}(b,a)-e^{-a-b}g_{i_1,i_2}(-a,-b)\bigg)
\end{align}
holds uniformly in complex numbers $\alpha,\beta\ll(\log\mathbf{q})^{-1}$, where $a=(\sigma_0+\alpha-1/2)\log\mathbf{q}$, $b=(\sigma_0+\beta-1/2)\log\mathbf{q}$. Moreover (\ref{3.6}) admits differentiations in $\alpha,\beta$.

\section{upper bound for $N_G(\mathcal{D},Q)$}
\label{sec4}
Let $\mathcal{D}_1$ be the closed rectangle with vertices $\sigma_0-iT$, $3-iT$, $\sigma_0+iT$, $3+iT$, where $\sigma_0=1/2-R/\log \mathbf{q}$.  Suppose $G(3+it,\chi)\neq0$. Determine $\text{arg}G(\sigma+iT,\chi)$ by continuation left from $3+iT$ and $\text{arg}G(\sigma-iT,\chi)$ by  continuation left from $3-iT$. If a zero is reached on the upper edge, use $\lim G(\sigma+iT+i\epsilon,\chi)$ as $\epsilon\rightarrow+0$ and $\lim G(\sigma-iT-i\epsilon,\chi)$ on the lower edge. Make horizontal cuts in $\mathcal{D}_1$ from the left side to the zeros of $G$ in $\mathcal{D}_1$. Applying the Littlewood's formula (to see \cite{Lit}), we have
\begin{align}
\label{2.41}
   &\int_{-T}^T\log|G(\sigma_0+it,\chi)|dt-\int_{-T}^T\log|G(3+it,\chi)|dt\notag\\
   &\ \ \ \ \ \ \ \ \ \ \ \ \ \ \ \ +\int_{\sigma_0}^3\text{arg}G(\sigma+iT,\chi)d\sigma-\int_{\sigma_0}^3\text{arg}G(\sigma-iT,\chi)d\sigma\notag\\
   &\ \ \ \ \ \ \ =2\pi\sum_{\rho\in \mathcal{D}_1}\text{dist}(\rho),
\end{align}
where $\text{dist}(\rho)$ is the distance of $\rho$ from the left side of $\mathcal{D}_1$.

Recall the definition of $\psi_i$ for $i=1,2$. A direct calculation shows that $\psi_i(s,\chi)\ll \mathbf{q}$ for $\text{Re}(s)>0$. Hence $G(s,\chi)\ll (T\mathbf{q})^2$ for $\text{Re}(s)>0$. Then we have
\begin{align}
\label{2.4.1}
   \int_{\sigma_0}^3\text{arg}G(\sigma+iT,\chi)d\sigma=O(\log (T\mathbf{q}))
\end{align}
by using Jensen's theorem in a familiar way as in \S9.4 of \cite{Tic}.
Similarly, we may have
\begin{align}
\label{2.4.2}
   \int_{\sigma_0}^3\text{arg}G(\sigma-iT,\chi)d\sigma=O(\log (T\mathbf{q})).
\end{align}
By a direct calculation we can see
\begin{align}
   \lambda L'(3+it,\chi)\psi_2(3+it,\chi)\ll O(1/\log\mathbf{q})\notag
\end{align}
and
\begin{align}
   |L(3+it,\chi)\psi_1(3+it,\chi)|\ge1-2\sum_2^\infty n^{-3}-\Big(\sum_2^\infty n^{-3}\Big)^2>1/3.\notag
\end{align}
Hence, from (\ref{2.1}) we have
\begin{align}
\label{2.51}
   \int_{-T}^T\log |G(3+it,\chi)|dt=\int_{-T}^T\log |L(3+it,\chi)\psi_1(3+it,\chi)|dt+O(T/\log\mathbf{q}).
\end{align}
Since for $\sigma>1$
\begin{align}
   \log L(s,\chi)=-\sum\frac{\Lambda(n)\chi(n)}{n^s\log n},\notag
\end{align}
it follows taking the real part that
\begin{align}
\label{2.611}
   \int_{-T}^T\log |L(3+it,\chi)|dt\ll 1.
\end{align}
For the entire function $\psi_1(s)$, it is easy to see, for $\sigma\ge3$,
\begin{align}
   |\psi_1(s)-1|\le\frac{1}{2^\sigma}+\frac1{3^\sigma}+\int_3^\infty\frac{\nu}{\nu^\sigma}\le\frac{1}{2^\sigma}+\frac52\frac{1}{3^\sigma}<2^{1-\sigma}.\notag
\end{align}
Therefore, $\log\psi_1(s)$ is analytic for $\sigma\ge3$. Integrating on the contour $\sigma+iT,~3\le\sigma<\infty$; $3+it,~-T\le t\le T;~\sigma+iT,~3\le \sigma<\infty$ gives
\begin{align}
\label{2.71}
   \int_{-T}^T\log|\psi_1(3+it,\chi)|dt\le\bigg|\int_{-T}^T\log\psi_1(3+it,\chi)dt\bigg|\le8\int_3^\infty\frac{d\sigma}{2^\sigma}=O(1).
\end{align}
Substituting (\ref{2.611}), (\ref{2.71}) into (\ref{2.51}), we have
\begin{align}
\label{2.4.3}
   \int_{-T}^T\log |G(3+it,\chi)|dt\ll T/\log\mathbf{q}.
\end{align}
From the above discussion, it is easy to see $G(3+it,\chi)\neq0$. Then using (\ref{2.4.1}), (\ref{2.4.2}), (\ref{2.4.3}) in (\ref{2.41}) and by the condition that $(\log Q)^6\le T\le(\log Q)^A$ for some constant $A\ge6$, we have
\begin{align}
   \int_{-T}^T\log|G(\sigma_0+it,\chi)|dt+O(T/\log \mathbf{q})=2\pi\sum_{\rho\in \mathcal{D}_1}\text{dist}(\rho).\notag
\end{align}
Since $\mathcal{D}\subset\mathcal{D}_1$ and all zeros of $G$ in closed rectangle $\mathcal{D}$ are at least distance $1/2-\sigma_0=R/\log\mathbf{q}$ from the line $\sigma=\sigma_0$, it follows that
\begin{align}
\label{2.11}
  N_G(\mathcal{D},\chi)\le\frac{\log\mathbf{q}}{2\pi R}\int_{-T}^T\log|G(\sigma_0+it,\chi)|dt+O(T).
\end{align}
It would be important if one can give a useful evaluation of the right side of (\ref{2.11}). However, at the present state of technology we are unable to obtain a useful evaluation. It is not until recently, due to Conrey, Iwaniec and Soundararajan's important work on Asymptotic Large Sieve, we can get a useful evaluation of this by averaging over the conductors $q$ and the primitive characters $\chi$.

Summing both sides of (\ref{2.11}) over $q\le Q$ and primitive characters $\chi$,  we have
\begin{align}
  N_G(\mathcal{D},Q)\le \frac{1}{4\pi R}\sum_{q\le Q}\frac{\Psi(q/Q)\log\mathbf{q}}{\phi(q)}\mathop{{\sum}^*}\limits_{\chi(\text{mod} q)}\int_{-T}^T\log|G(\sigma_0+it,\chi)|^2dt+o(N(T,Q)).\notag
\end{align}
By the concavity of the Logarithm function we get
\begin{align}
\label{2.5}
   N_G(\mathcal{D},Q)\le\Big(\frac1{2R}\log c(\theta,r,R)+o(1)\Big) N(T,Q),
\end{align}
where $c(\theta,r,R)$ is defined by
\begin{align}
\label{2.7}
   c(\theta,r,R)\sum_{q\le Q}\Psi(\frac qQ)\log\mathbf{q}\frac{\phi^*(q)}{\phi(q)}=\frac1{2T}\sum_{q\le Q}\frac{\Psi(q/Q)\log\mathbf{q}}{\phi(q)}\mathop{{\sum}^*}\limits_{\chi(\text{mod} q)}\int_{-T}^T|G(\sigma_0+it,\chi)|^2dt.
\end{align}
Here $\phi^*(q)$ denotes the number of primitive characters (mod $q$). In the process to obtain (\ref{2.5}), we have used $\log \mathbf{q}T$ to replace $\log\mathbf{q}$. We can make this replacement because we have assumed that $(\log Q)^6\le T\le(\log Q)^A$ for some constant $A\ge6$.

It is obvious that $\Psi(\frac qQ)\log\mathbf{q}$ in (\ref{2.7}) still satisfies the condition of  $\Psi(\frac qQ)$ given in the definition of $N(T, Q)$. Hence we can evaluate the right side of (\ref{2.7}) by using the Asymptotic Large Sieve result. Recalling the definition of $G(s,\chi)$ in (\ref{2.1}), $\lambda=1/(r\log\mathbf{q})$, we have
\begin{align}
   |G(\sigma_0+it,\chi)|^2&=L(\sigma_0+it,\chi)L(\sigma_0-it,\bar{\chi}) \psi_{1}(\sigma_0+it,\chi)\psi_{1}(\sigma_0-it,\bar{\chi})\notag\\
   &+\frac1{r\log\mathbf{q}}L'(\sigma_0+it,\chi)L(\sigma_0-it,\bar{\chi}) \psi_{2}(\sigma_0+it,\chi)\psi_{1}(\sigma_0-it,\bar{\chi})\notag\\
   &+\frac1{r\log\mathbf{q}} L(\sigma_0+it,\chi)L'(\sigma_0-it,\bar{\chi}) \psi_{1}(\sigma_0+it,\chi)\psi_{2}(\sigma_0-it,\bar{\chi})\notag\\
   &+\frac1{r^2\log^2\mathbf{q}} L'(\sigma_0+it,\chi)L'(\sigma_0-it,\bar{\chi}) \psi_{2}(\sigma_0+it,\chi)\psi_{2}(\sigma_0-it,\bar{\chi}).\notag
\end{align}
Substituting this into (\ref{2.7}), we have
\begin{align}
\label{h3.1}
   c(\theta,r,R)\sum_{q\le Q}\Psi(\frac qQ)\log\mathbf{q}\frac{\phi^*(q)}{\phi(q)}=\frac1{2T}\Bigg(D_{11}+\frac1{r\log\mathbf{q}}D_{21} +\frac1{r\log\mathbf{q}}D_{12}+\frac1{r^2\log^2\mathbf{q}}D_{22}\Bigg),
\end{align}
where
\begin{align}
   &D_{i_1i_2}=\sum_{q\le Q}\frac{\Psi(q/Q)\log\mathbf{q}}{\phi(q)}\notag\\
   &\ \ \ \ \ \ \ \ \ \ \ \ \ \times\mathop{{\sum}^*}\limits_{\chi(\text{mod} q)}\int_{-T}^TL^{(i_1-1)}(\sigma_0+it,\chi)L^{(i_2-1)}(\sigma_0-it,\bar{\chi}) \psi_{i_1}(\sigma_0+it,\chi)\psi_{i_2}(\sigma_0-it,\bar{\chi})dt.\notag
\end{align}
Choosing $\alpha=\beta=0$ in (\ref{3.6}), we get
\begin{align}
\label{+1}
   D_{11} &\sim 2T \bigg(\sum_{q\le Q}\Psi(\frac qQ)\log\mathbf{q}\frac{\phi^*(q)}{\phi(q)}\bigg) \bigg(\frac{g_{1,1}(b,a)-e^{-a-b}g_{1,1}(-a,-b)}{(a+b)\theta}\bigg)\bigg|_{a=b=-R}
\end{align}
with $g_{i_1,i_2}(a,b)$ defined by (\ref{3.5}).
Differentiating (\ref{3.6}) in $\alpha$ and choosing $\alpha=\beta=0$, we get
\begin{align}
\label{+2}
   D_{21} &\sim 2T\log\mathbf{q} \bigg(\sum_{q\le Q}\Psi(\frac qQ)\log\mathbf{q}\frac{\phi^*(q)}{\phi(q)}\bigg) \partial_a\bigg(\frac{g_{2,1}(b,a)-e^{-a-b}g_{2,1}(-a,-b)}{(a+b)\theta}\bigg)\bigg|_{a=b=-R}.
\end{align}
Differentiating (\ref{3.6}) in $\beta$ and choosing $\alpha=\beta=0$, we get
\begin{align}
\label{+3}
   D_{12} &\sim 2T\log\mathbf{q} \bigg(\sum_{q\le Q}\Psi(\frac qQ)\log\mathbf{q}\frac{\phi^*(q)}{\phi(q)}\bigg) \partial_b\bigg(\frac{g_{1,2}(b,a)-e^{-a-b}g_{1,2}(-a,-b)}{(a+b)\theta}\bigg)\bigg|_{a=b=-R}.
\end{align}
Differentiating (\ref{3.6}) both in $\alpha,~\beta$ and choosing $\alpha=\beta=0$, we get
\begin{align}
\label{+4}
   D_{22} &\sim 2T\log^2\mathbf{q} \bigg(\sum_{q\le Q}\Psi(\frac qQ)\log\mathbf{q}\frac{\phi^*(q)}{\phi(q)}\bigg) \partial_{a,b}\bigg(\frac{g_{2,2}(b,a)-e^{-a-b}g_{2,2}(-a,-b)}{(a+b)\theta}\bigg)\bigg|_{a=b=-R}.
\end{align}
Substituting (\ref{+1})-(\ref{+4}) into (\ref{h3.1}), we have
\begin{align}
   c(\theta,r,R)\sim&h_{1,1}(a,b)\big|_{a=b=-R}+\frac1r\partial_ah_{2,1}(a,b)\big|_{a=b=-R}\notag\\
     &+\frac1r\partial_bh_{1,2}(a,b)\big|_{a=b=-R}+\frac1{r^2}\partial_{a,b}h_{2,2}(a,b)\big|_{a=b=-R},
\end{align}
where
\begin{align}
   h_{i_1,i_2}(a,b)=&\frac1{\theta(a+b)}\Bigg\{\int_0^1\Big(P_{i_1}'(t)+b\theta P_{i_1}(t)\Big)\Big(P_{i_2}'(t)+a\theta P_{i_2}(t)\Big)dt\notag\\
   &-e^{-a-b}\int_0^1\Big(P_{i_1}'(t)-a\theta P_{i_1}(t)\Big)\Big(P_{i_2}'(t)-b\theta P_{i_2}(t)\Big)dt\Bigg\}.
\end{align}
Taking $\theta=1-\epsilon$, $r=1.154$, $R=0.617$,
\begin{align}
   &P_1(x)=x-0.158x(1-x)+0.25x^2(1-x),\notag\\
   &P_2(x)=x-0.492x(1-x)+0.075x^2(1-x),\notag
\end{align}
and making $\epsilon\rightarrow0$, we get $c(\theta,r,R)=1.230108\cdots$, then from (\ref{2.5}) we have
\begin{align}
   N_G(\mathcal{D},Q)\le 0.167835N(T,Q).\notag
\end{align}

\section{zeros of $\xi'$ on the critical line}
\label{sec5}
In this section, we will study the number of zeros of the family of $\xi'(s,\chi)$ on the critical line by using Asymptotic Large Sieve and Levinson's method. Actually, we do not follow the way in \cite{Con-Zero}, which studied the proportions of zeros of $\xi^{(n)}(s,1)$, $n\ge0$, on the critical line, but use the Levinson's method generalized by Conrey in \cite{Con-More}. One may see that the way in this section will give a better result about the proportion of zeros of the family of $\xi'(s,\chi)$ on the critical line than the way in \cite{Con-Zero}.

From the functional equation (\ref{s2.3}) it is easy to see that $h(\chi)\xi^{(n)}(s,\chi)$ is real for $s=1/2+it$ when $n$ is even and is purely imaginary when $n$ is odd. Let $\delta\neq0$ be real, $g_n, ~n\ge1$, be complex numbers with $g_n$ real if $n$ is even and $g_n$ purely imaginary if $n$ is odd. Let $T$ be a large parameter and
\begin{align}
   \mathcal{L}=\log (T\mathbf{q}).\notag
\end{align}
Now define
\begin{align}
   \eta(s,\chi)=(1-\delta)\xi(s,\chi)+\delta\xi'(s,\chi)\mathcal{L}^{-1}+\sum_{n=1}^Ng_n\xi^{(n)}(s,\chi)\mathcal{L}^{-n}\notag
\end{align}
for some fixed $N$. Then, for $s=1/2+it$,
\begin{align}
   \delta h(\chi)\xi'(s,\chi)=\text{Im}\Big(h(\chi)\eta(s,\chi)\Big),\notag
\end{align}
so that $\xi'(s)=0$ on $\sigma=1/2$ if and only if $\text{Im}\Big(h(\chi)\eta(s,\chi)\Big)=0$. Observe that for every change of $\pi$ in the argument of $h(\chi)\eta(s,\chi)$ it must be the case that $\text{Im}\Big(h(\chi)\eta(s,\chi)\Big)$ has at least one zero. Hence it follows that
\begin{align}
\label{s3.3}
   N_{\xi',c}(T,\chi)\ge\frac1\pi\Delta_C\text{arg}\Big(h(\chi)\eta(s,\chi)\Big)=\frac1\pi\Delta_C\text{arg}\eta(s,\chi),
\end{align}
where $\Delta_C\text{arg}$ stands for the variation of the argument as $s$ runs over the critical line from $1/2-iT$ to $1/2+iT$ passing the zeros of $\eta(s,\chi)$ from the east side.

To estimate the change in argument of $\eta(s,\chi)$ on the critical line, we let $\eta(s,\chi)=H(s,\chi)V(s,\chi)$, where $H(s,\chi)$ is defined in (\ref{s2.2}) and
\begin{align}
   V(s,\chi)=&(1-\delta)L(s,\chi)+\frac{\delta}{\mathcal{L}}\bigg(\frac{H'}{H}(s,\chi)L(s,\chi)+ L'(s,\chi)\bigg)\notag\\
   &+\sum_{n=1}^N\frac{g_n}{\mathcal{L}^{n}}\sum_{k=0}^n\binom{n}{k}\frac{H^{(n-k)}(s,\chi)}{H(s,\chi)}L^{(k)}(s,\chi). \notag
\end{align}
By the Stirling formula, for $|t|\ge2$, we have
\begin{align}
   \text{arg}H(1/2+it)=\frac t2\log\frac{|qt|}{2\pi e}+O(1),\notag
\end{align}
and
\begin{align}
   \frac{H^{(m)}}{H}(s,\chi)=(1/2\log\frac {qs}{2\pi})^m(1+O(1/|t|))\notag
\end{align}
for $t\ge10$, $0<\sigma<A_1$, here $A_1$ can be any positive constant. (For a proof of these formulas, see Lemma 1 of \cite{Con-Zero}.) Hence we may have
\begin{align}
\label{s3.2}
   \Delta\text{arg}\eta(1/2+it,\chi)\big|_{-T}^T=T\log\frac{qT}{2\pi e}+\Delta\text{arg}V(1/2+it,\chi)\big|_{-T}^T+O(T)
\end{align}
and denote $V(s,\chi)$ by
\begin{align}
   V(s,\chi)=\Bigg\{\Bigg(1-\delta+\delta\bigg(\frac{\log \frac{qs}{2\pi}}{2\mathcal{L}}+\frac1{\mathcal{L}}\frac{d}{ds}\bigg)Q_0\bigg(\frac{\log \frac{qs}{2\pi}}{2\mathcal{L}}+\frac1{\mathcal{L}}\frac{d}{ds}\bigg)\Bigg)L(s,\chi)\Bigg\}(1+O(1/|t|))\notag
\end{align}
with
\begin{align}
   Q_0(x)=1+\sum_{n=1}^N\frac{g_n}{\delta}x^{n-1}.\notag
\end{align}
As in (\ref{2.30}), we use the mollifier
\begin{align}
   \psi(s,\chi)=\sum_{n\le X}\frac{\mu(n)\chi(n)}{n^{s+R/\mathcal{L}}}P\bigg(\frac{\log X/n}{\log X}\bigg),\notag
\end{align}
where $\mu$ is the M\"{o}bius function, $X=\mathbf{q}^\theta$ with $0<\theta<1$. Here $P$ is a polynomial with $P(0)=0$, $P(1)=1$ which will be specified later.
By the Cauchy's argument principle and using Jenson's theorem as in section \ref{sec4}, we may have
\begin{align}
    \Big|\Delta\text{arg}V(1/2+it,\chi)\big|_{-T}^T\Big|=2\pi N_{\xi'}(\mathcal{D},\chi)+O(\mathcal{L}).\notag
\end{align}
Here $N_{\xi'}(\mathcal{D},\chi)$ denote the number of $\xi'(s,\chi)$ in the closed rectangle $\mathcal{D}$ defined before.
Then, if $Q_0(1/2)=2$, by applying Jenson's theorem and Littlewood's formula as in section \ref{sec4}, we can show that
\begin{align}
\label{s3.1}
   \Big|\Delta\text{arg}V(1/2+it,\chi)\big|_{-T}^T\Big|\le \frac{\log\mathbf{q}}{R}\int_{-T}^T\log\big|V\psi(\sigma_0+it,\chi)\big|dt+O(T),
\end{align}
where
\begin{align}
   \sigma_0=1/2-R/\log \mathbf{q},\notag
\end{align}
$R\ll1$ is a positive real number. Here, the reason for requiring the condition $Q_0(1/2)=2$ is to ensure the integration $\int_{-T}^T\log |V\psi(3+it,\chi)|dt$ caused by using Littlewood's formula is $O(T/\log\mathbf{q})$. To evaluate that the integration in the right of (\ref{s3.1}) we use the following useful approximation to $V(s,\chi)$,
\begin{align}
\label{s3.4}
   U(s,\chi)=\bigg(1-\delta+\delta\Big(1+\frac2{\mathcal{L}}\frac{d}{ds}\Big)Q\Big(-\frac1{\mathcal{L}}\frac{d}{ds}\Big)\bigg)L(s,\chi),
\end{align}
where
\begin{align}
   Q(x)=\frac12Q_0(1/2-x).\notag
\end{align}
It is easy to see that the error caused by the substitution of $V$ with $U$ can be absorbed into the error term in (\ref{s3.1}). If restrict $Q(x)$ to be a real polynomial, we can see that the restriction of $g_n$ and the condition that $Q_0(1/2)=2$ are equivalent to $Q'(x)=Q'(1-x)$ and $Q(0)=1$. Hence by Lemma \ref{lems1} and (\ref{s3.3})-(\ref{s3.4}),
\begin{align}
   N_{\xi',c}(T,\chi)\ge N(T,\chi)(1+o(1))-\frac{\log\mathbf{q}}{\pi R}\int_{-T}^T\log|U\psi(\sigma_0+it,\chi)|dt.\notag
\end{align}
Summing both sides of the above formula over $q\le Q$ and primitive characters $\chi$,  we have
\begin{align}
\label{s3.5}
  N_{\xi',c}(T,Q)\ge N(T,Q)(1+o(1))-\frac{1}{2\pi R}\sum_{q\le Q}\frac{\Psi(\frac qQ)\log \mathbf{q}}{\phi(q)}\mathop{{\sum}^*}\limits_{\chi(\text{mod} q)}\int_{-T}^T\log|U\psi(\sigma_0+it,\chi)|^2dt.
\end{align}

Recall that $(\log Q)^6\le T\le(\log Q)^A$ for some constant $A\ge6$. We refine the definition of $U$ which was previously defined in (\ref{s3.4}) to
\begin{align}
\label{s3.6}
   U(s,\chi)=\bigg(1-\delta+\delta\Big(1+\frac2{\log\mathbf{q}}\frac{d}{ds}\Big)\bigg)Q\bigg(-\frac1{\log\mathbf{q}}\frac{d}{ds}\bigg)L(s,\chi)
\end{align}
with $Q$ defined as before. Then, in (\ref{s3.5}), the error caused by this refinement can be absorbed by the error term.
Substituting this refined $U$ into (\ref{s3.5}) and using the concavity of the Logarithm function, we get
\begin{align}
\label{s3.7}
    N_{\xi',c}(T,Q)\ge\Big(1-\frac1R\log c_1(\theta,R)+o(1)\Big) N(T,Q),
\end{align}
where $c_1(\theta,R)$ satisfies
\begin{align}
\label{s3.8}
   c_1(\theta,R)\sum_{q\le Q}\Psi(\frac qQ)\log\mathbf{q}\frac{\phi^*(q)}{\phi(q)}=\frac1{2T}\sum_{q\le Q}\frac{\Psi(\frac qQ)\log\mathbf{q}}{\phi(q)}\mathop{{\sum}^*}\limits_{\chi(\text{mod} q)}\int_{-T}^T|U\psi(\sigma_0+it,\chi)|^2dt.
\end{align}

Then, using formula (\ref{3.6}) in (\ref{s3.8}) as in section \ref{sec4}, we get
\begin{align}
\label{s3.11}
   c_1(\theta,R)\sim&\bigg(1-\delta+\delta(1+2\partial_a)Q(-\partial_a)\bigg)\bigg(1-\delta+\delta(1+2\partial_b)Q(-\partial_b)\bigg)\notag\\
   &\times\Bigg(\frac{g(b,a)-e^{-a-b}g(-a,-b)}{\theta(a+b)}\Bigg)\Bigg|_{a=b=-R},
\end{align}
where
\begin{align}
   g(a,b)=\int_0^1\Big(P'(t)+a\theta P(t)\Big)\Big(P'(t)+b\theta P(t)\Big)dt.\notag
\end{align}
Taking $\theta=1-\epsilon$, $R=0.746$, $\delta=0.771$,
\begin{align}
   P(x)= x-0.482x(1-x)-0.392x^2(1-x)-0.262x^3(1-x),\notag
\end{align}
\begin{align}
   Q(x)=1-0.673x+0.369(x^2/2-x^3/3)-4.635(x^3/3-x^4/2+x^5/5)\notag
\end{align}
into (\ref{s3.11}) and making $\epsilon\rightarrow0$, we have by (\ref{s3.8}) that
\begin{align}
    N_{\xi',c}(T,Q)\ge0.93828 N(T,Q).\notag
\end{align}

By now we have proved that more than $93.828\%$ of zeros of the family of $\xi'(s,\chi)$ are on the critical line, which gives a lower bound for $N_{\xi',c}(T,Q)$. We note that the result obtained by the method in \cite{Con-Zero} is actually equal to the case $\delta=1$ here. In \cite{Con-Zero}, it proved that more than 81.37\% of zeros of $\xi'(s,1)$ are on the critical line with a mollifier of length $T^{1/2-\epsilon}$. To obtain a lower bound for $N_{\xi',c}(T,Q)$, we can also follow the way in \cite{Con-Zero} with a mollifier of length $\mathbf{q}^{1-\epsilon}$ as used in this section, and then the result we can obtain is $N_{\xi',c}(T,Q)\ge 0.8429N(T,Q)$.

\section{completion of the proof}

We have obtained that
\begin{align}
  N_G(\mathcal{D},Q)\le 0.167835N(T,Q)\notag
\end{align}
in section \ref{sec4} and
\begin{align}
    N_{\xi',c}(T,Q)\ge0.93828 N(T,Q).\notag
\end{align}
in section \ref{sec5}.
Then by (\ref{2.2}), (\ref{2.21}) we have
\begin{align}
   N_d(T,Q)&\ge \frac12N(T,Q)+\frac12N_{\xi',c}(T,Q)-N_G(\mathcal{D},Q)\notag\\
   &\ge(\frac12+0.46914-0.167835)N(T,Q)>0.8013N(T,Q),\notag
\end{align}
and for simple zero
\begin{align}
   N_s(T,Q)&\ge N_{\xi',c}(T,Q)-2N_G(\mathcal{D},Q)\notag\\
   &\ge(0.93828-0.33567)N(T,Q)=0.60261N(T,Q).\notag
\end{align}
Hence we have proved Theorem \ref{thm1}.

If assume the Generalized Riemann Hypothesis holds, from (\ref{2.3}), (\ref{2.31}) we have
\begin{align}
   N_d(T,Q)\ge N(T,Q)-N_G(\mathcal{D},Q)\ge 0.83216N(T,Q),\notag
\end{align}
and
\begin{align}
   N_s(T,Q)\ge N(T,Q)-2N_G(\mathcal{D},Q)\ge 0.66433N(T,Q).\notag
\end{align}
Hence we have proved Theorem \ref{thm2}.

{\bf Remark.} Using Levinson's method, we may estimate the proportions for simple zeros of the family of
$\xi^{(n)}(s,\chi)$, $n\ge0$, and then use Farmer's combination method as in \cite{Far} to get a lower bound for the proportion of distinct zeros only. However, the result obtained by this way is much worse than Theorem \ref{thm1}. As mentioned in \cite{CIS-Crit}, our method can be generalized to $GL_2$ and $GL_3$ $L$-function.

\end{document}